\documentclass[12pt]{amsart}
\usepackage{amscd,amssymb}
\usepackage[graph,frame,poly,arc]{xy}  
\usepackage[plainpages,backref,urlcolor=blue]{hyperref}

\theoremstyle{plain}
\newtheorem{thm}[subsection]{Theorem}
\newtheorem{lem}[subsection]{Lemma}
\newtheorem{prop}[subsection]{Proposition}
\newtheorem{cor}[subsection]{Corollary}

\theoremstyle{definition}
\newtheorem{rk}[subsection]{Remark}

\newtheorem{ex}[subsection]{Example}

\numberwithin{equation}{section}
\newcommand{\OO}{{\mathcal O}}

\newcommand{\I}{{\mathcal I}}

\newcommand{\F}{{\mathcal F}}
\newcommand{\A}{{\mathcal A}}
\newcommand{\B}{{\mathcal B}}

\newcommand{\R}{\mathbb{R}}
\newcommand{\C}{\mathbb{C}}
\newcommand{\K}{\mathbb{K}}
\newcommand{\PP}{\mathbb{P}}

\begin{document}
%\date{June 4, 2009}

\title [On supersolvable and nearly supersolvable line arrangements]
{On supersolvable and nearly supersolvable line arrangements}

\author[Alexandru Dimca]{Alexandru Dimca$^1$}
\address{Universit\'e C\^ ote d'Azur, CNRS, LJAD, France }
\email{dimca@unice.fr}

\author[Gabriel Sticlaru]{Gabriel Sticlaru}
\address{Faculty of Mathematics and Informatics,
Ovidius University
Bd. Mamaia 124, 900527 Constanta,
Romania}
\email{gabrielsticlaru@yahoo.com }

%\thanks{$^1$ Partially supported by Institut Universitaire de France.} 

\thanks{$^1$ This work has been supported by the French government, through the $\rm UCA^{\rm JEDI}$ Investments in the Future project managed by the National Research Agency (ANR) with the reference number ANR-15-IDEX-01.}



\begin{abstract} We introduce a new class of line arrangements in the projective plane, called nearly supersolvable, and show that any arrangement in this class is either free or nearly free. More precisely, we show that the minimal degree of a Jacobian syzygy for the defining equation of the line arrangement, which is a subtle algebraic invariant, is determined in this case by the combinatorics. When such a line arrangement is nearly free, we discuss the splitting types and the jumping lines of the associated rank two vector bundle, as well as the corresponding jumping points, introduced recently by S. Marchesi and J. Vall\` es. As a by-product of our results, we get a version of the Slope Problem, valid over the real and the complex numbers as well.

\end{abstract}
 
\maketitle

\section{Introduction} 
Let $\A:f=0$ be a line arrangement  in the complex projective plane $\PP^2$.
An intersection point $p$ of $\A$ is called a {\it modular} point if for any other intersection point $q$ of $\A$, the line $\overline {pq}$ determined by the points $p$ and $q$ belongs to the arrangement $\A$. The arrangement $\A$ is {\it supersolvable} if it has a modular intersection point.
Supersolvable arrangements have many interesting properties, in particular they are free arrangements, see \cite{DHA,T1,T2} or  \cite[Prop 5.114]{OT} and \cite[Theorem 4.2]{JT}).

In this note we introduce a new class of line arrangements as follows.  An intersection point $p$ of a line arrangement $\A$ in $\PP^2$ is called a {\it nearly modular} point of $\A$ if the following two properties hold.
\begin{enumerate}

\item For any  intersection point $q\ne p$ of $\A$, with the exception of a unique double point $p' \ne p$ of $\A$, the line $\overline {pq}$ determined by the points $p$ and $q\ne p'$ belongs to the arrangement $\A$.

\item The line $L=\overline {pp'}$ is not in $\A$ and contains only two multiple points  of $\A$, namely $p$ and $p'$.

\end{enumerate}

 The arrangement $\A$ is {\it nearly supersolvable} if $\A$ is not supersolvable, but it has a nearly modular intersection point $p$. For any pair $(\A,p)$, with $\A$ a nearly supersolvable arrangement and $p$ a nearly modular point of $\A$, we get a supersolvable line arrangement $\B=\B(\A,p)$, by adding to $\A$ the line $L=\overline {pp'}$. 
 
 In the second section we recall the definition of the minimal degree $mdr(f)$ of a Jacobian syzygy for $f$, as well as the definition and some basic properties of the free and nearly free line arrangements. The only new result here is Proposition \ref{propNF} which gives a new view point on the jumping point of a nearly free arrangement, a notion introduced by 
S. Marchesi and J. Vall\` es in \cite{MaVa}. In fact our result was motivated and inspired by \cite[Theorem 2.1]{MaVa}, see Remark \ref{rkNF} for more details on the relation between these two results.

In the third section we obtain the relation between the minimal degree $mdr(f)$ and the multiplicity of a modular point of $\A$ in Proposition \ref{propA} and we introduce a number of line arrangements to illustrate our results.
In the fourth section we prove the main result, Theorem \ref{thmB}, saying that the multiplicity of a nearly modular point determines the minimal degree $mdr(f)$ as well as whether the nearly supersolvable arrangement $\A$ is free or nearly free. Example \ref{exB}, Remark \ref{rkNF3} and Example \ref{exC} illustrate 
Theorem \ref{thmB} and the properties of the jumping points of nearly free arrangements. As a by-product we get the following version of the Slope Problem, valid over $\K=\R$ and $\K=\C$ as well.

\begin{thm}
\label{thmSlope}
For a configuration of $n$ points in the affine plane $\K^2$, not all of them on the same line and such that there exist two of them, say $P_1$ and $P_2$,  determining a line of unique slope, then the number of distinct slopes of the lines determined by the $n$ points is at least $n$.
\end{thm}
For more on the Slope Problem, and a precise statement of Theorem \ref{thmSlope} we refer to \cite{Scott,T1,U} and Theorem \ref{thmSlope2} below.
In the final section we consider  the sheaf $T\langle \A \rangle $ of logarithmic vector fields along the nearly supersolvable line arrangement $\A$ and investigate its splitting types and its jumping lines, using  a key result due to S. Marchesi and J. Vall\` es in \cite{MaVa}.

\bigskip
We would like to thank the referees for their very useful remarks which
helped us to improve both the presentation and the results in our manuscript.

\section{Free and nearly free line arrangements} 

Let $S=\C[x,y,z]$ be the polynomial ring in three variables $x,y,z$ with complex coefficients, and let $\A:f=0$ be an arrangement of  $d$ lines in the complex projective plane $\PP^2$. The minimal degree of a Jacobian syzygy for the polynomial $f$ is the integer $mdr(f)$
defined to be the smallest integer $m\geq 0$ such that there is a nontrivial relation
\begin{equation}
\label{rel_m}
 af_x+bf_y+cf_z=0
\end{equation}
among the partial derivatives $f_x, f_y$ and $f_z$ of $f$ with coefficients $a,b,c$ in $S_m$, the vector space of  homogeneous polynomials in $S$ of degree $m$. When $mdr(f)=0$, then $\A$ is a union of $d$ lines passing through one point, a situation easy to analyse. We assume from now on in this note that 
$$ mdr(f)\geq 1.$$
It was shown by Ziegler \cite{Zi}, see also for details \cite[Remark 8.5]{DHA}, that this algebraic  invariant $mdr(f)$ is not determined by the combinatorics of the line arrangement $\A:f=0$ in general.
Denote by $\tau(\A)$ the global Tjurina number of the arrangement $\A$, which is the sum of the Tjurina numbers $\tau(\A,a)$ of the singular points $a$ of $\A$. If $n_k$ is the number of intersection points in $\A$ of multiplicity $k$, for $k \geq 2$, then one has
\begin{equation} \label{r0}
\tau(\A)= \sum_{k \geq 2}n_k(k-1)^2.
\end{equation}
Indeed, any singular point $a$ of multiplicity $k \geq 2$ of a line arrangement $\A$ being weighted homogeneous, the local Tjurina number $\tau(\A,a)$ coincides with the local Milnor number $\mu(\A,a)=(k-1)^2$, see \cite{KS}.
Moreover, one has  
\begin{equation} \label{r1}
\tau(\A) \leq (d-1)^2-r(d-r-1),
\end{equation}
where $r=mdr(f)$, see \cite{Dmax,duPCTC}, and equality holds if and only if the line arrangement  $\A:f=0$ is free. 
In this case, $d_1=r$ and $d_2=d-1-r$ are the exponents of the free arrangement $\A$. Note that for any free line arrangement one has $d_1=r \leq d-1-r=d_2$, and hence $r<d/2$ in this case. Usually, the free arrangements are defined as follows. Let $AR(f) \subset S^3$ be the graded $S$-module such, for any integer $m$, the corresponding homogeneous component $AR(f)_m$ consists of all the triples $\rho=(a,b,c)\in S^3_m$ satisfying \eqref{rel_m}. Then the arrangement $\A:f=0$ is said to be free if the graded $S$-module $AR(f)$ is free. In such a situation, one has $AR(f)=S(-d_1)\oplus S(-d_2)$, where $(d_1,d_2)$ are the exponents of $\A$ as defined above. The associated coherent sheaf on $\PP^2$ of the graded module $AR(f)$ is just  $T\langle \A \rangle(-1) $,  where $T\langle \A \rangle $ is the sheaf of logarithmic vector fields along $\A$ as considered for instance in \cite{AD, DS14}. For a free arrangement $\A$ as above, this yields
$$T\langle \A \rangle(-1) =\OO_{\PP^2}(-d_1)\oplus \OO_{\PP^2}(-d_2).$$
For basic facts on free arrangements, please refer to \cite{DHA,OT, Yo}.

Similarly, the line arrangement  $\A:f=0$ is nearly free, a notion introduced in \cite{DStRIMS} motivated by the study of rational cuspidal curves in \cite{DStFD}, if and only if 
\begin{equation} \label{r2}
\tau(\A)=(d-1)^2-r(d-r-1)-1,
\end{equation}
where $r=mdr(f)$, see \cite{Dmax}. In this case, $d_1=r$ and $d_2=d-r$ are the exponents of the nearly free arrangement $\A$, and $d_1=r \leq d-r=d_2$. Therefore $2r \leq d$ in this case. In terms of the graded module $AR(f)$, the nearly free arrangements are described by the following result.

\begin{prop}
\label{propNF}
Let $\A:f=0$  be
an arrangement  of $d$ lines in $\PP^2$ and let $r=mdr(f)$. Then, for any choice of a nonzero syzygy
$\rho_1 \in AR(f)_r$, there is a homogeneous ideal $I \subset S$ and an exact sequence
$$0 \to S(-r) \to AR(f) \to I(r-d+1) \to 0,$$
such that the following hold.
\begin{enumerate}

\item The ideal $I$ is saturated, defines a subscheme of $\PP^2$ of dimension at most 0, and its degree is given by
$$\deg I=(d-1)^2-r(d-r-1)-\tau(\A).$$

\item The line arrangement $\A$ is free if and only if $I=S$.

\item The line arrangement $\A$  is nearly free if and only if $I$ defines a reduced point $P(\A)$ in $\PP^2$. The exact sequence and the point $P(\A)$ are unique when $2r<d$, i.e. when the exponents of the nearly free arrangement $\A$ satisfy $r=d_1<d_2=d-r$.

\end{enumerate}
\end{prop} 
\proof
The proof follows from the exact sequence (3.3) in \cite{Dmax}, if we define $I(r-d+1)$ to be the image of the morphism $v$. The claim that $I$ is saturated follows from the fact that the graded $S$-module $AR(f)$ is clearly saturated, and by using the long exact sequence of cohomology groups coming from the exact sequence \eqref{exseq} below and the vanishing
$H^1(\PP^2, \OO_{\PP^2}(m))=0$ for any integer $m$. The claim about the degree $\deg I$ follows from the equality 
$$\deg I= \dim S_m/I_m$$
for $m>>0$, and a direct computation of $\dim S_m/I_m$ using the the exact sequence of graded $S$-modules above and the obvious exact sequence
$$ 0 \to AR(f) \to S^3 \to J_f (d-1)\to 0$$
where $J_f$ is the Jacobian ideal of $f$, i.e. the ideal generated by $f_x,f_y,f_z$ in $S$.
More precisely, if $M(f)=S/J_f$ denotes the corresponding Jacobian algebra, one has
\begin{equation} \label{degI}
\dim S_m/I_m=(d-1)^2-r(d-r-1)-\dim M(f)_{m+2(d-1)-r},
\end{equation}
for any $m \geq \max\{0, 2r-1-d\}$. The claim follows, since $\tau(\A)=\dim M(f)_s$ for $s>>0$.
The claim that the ideal $I$ defines a simple point on $\PP^2$ if and only if $\A$ is nearly free follows from \cite[Theorem 4.1]{Dmax}. Note that the last equality in this result has a minor misprint, the correct version is $\delta(f)_{d-r}=2$.
\endproof

\begin{rk}
\label{rkNF}
Note that Proposition \ref{propNF} holds for any reduced plane curve $C:f=0$ with exactly the same proof, see \cite{DStJump} for further results in this general setting. If we consider the associated coherent sheaves, the exact sequence in 
Proposition \ref{propNF} becomes
\begin{equation} \label{exseq}
0 \to \OO_{\PP^2}(-r) \to T\langle \A \rangle (-1) \to \I(r-d+1) \to 0.
\end{equation}
This exact sequence appeared first in \cite[Theorem 2.1]{MaVa} in the case of a nearly free curve, and this result was the motivation and the inspiration for our approach. It is known that for any finitely generated graded $S$-module
$F$, there is an associated coherent sheaf $\F$ on $\PP^2$, and conversely, for any coherent sheaf $\F$ on $\PP^2$, the graded $S$-module
$$\Gamma(\F)= \sum_{k\geq 0}H^0(\PP^2,\F(k))$$
is finitely generated. However, the transformation
$$F \to \F \to \Gamma(\F)$$
is not the identity, i.e. the graded module $F$ cannot be recovered from the associated coherent sheaf $\F$, though one has $F_k=H^0(\PP^2,\F(k))$ for $k$ large enough. Due to this fact, it seems to us that a statement about graded modules is not just a translation of a statement about their associated coherent sheaves, but it is slightly more precise.

Following \cite{MaVa}, we call $P(\A)$ the jumping point of the nearly free arrangement $\A$. S. Marchesi and J. Vall\`es in \cite{MaVa} have considered this jumping point only when $2r<d$, 
and in this case $P(\A)$ is determined by $\A$. In fact, when $2r=d$, the corresponding vector bundle $T\langle \A \rangle$ is a twist of the tangent bundle of $\PP^2$, and hence it has no jumping lines.
We prefer to consider the jumping point $P(\A)$ even in the case $2r=d$, in spite of the fact  that it does not create any jumping line.
The general situation of a reduced non free plane curve $C$ is discussed in \cite{DStJump}, where the jumping point is replaced by a jumping 0-dimensional subscheme of $\PP^2$, whose relation with the jumping lines of the corresponding vector bundle $T\langle C \rangle$ is rather subtle.
\end{rk}

\begin{rk}
\label{rkNF0}
When $\A:f=0$ is nearly free with exponents $(d_1,d_2)$, then one has the following explicit description of the ideal $I$, which occurs already in \cite{MaVa}. Let $\rho_i \in AR(f)_{d_i}$ for $i=1,2,3$ a minimal set of generators for the graded $S$-module $AR(f)$, where $d_3=d_2$. Then there is a relation
$$h_1\rho_1+h_2\rho_2+h_3\rho_3=0$$
where $h_1 \in S$ is homogeneous of degree $d_2-d_1+1$ and $h_2,h_3$ are linearly independent linear forms in $S$, see \cite{DStRIMS}. With this notation, the ideal $I$ is generated by $h_2$ and $h_3$, see \cite{Dmax}, the discussion following the equation (3.3).
Hence, when $d_1<d_2$, the jumping point $P(\A)$ is defined by the system of equations $h_2=h_3=0$.
For a concrete situation, see Example \ref{exB} below.
\end{rk}

\section{Multiplicity of modular points and minimal degree of Jacobian syzygies}

Recall the following result, see \cite[Lemma2.1]{T2}. We denote by $m_p(\A)$ the multiplicity of an intersection point $p$ of $\A$, that is the number of lines in $\A$ passing through the point $p$.

\begin{lem}
\label{lemA}
If $\A$ is a supersolvable line arrangement, $p$ a modular point of $\A$, and $q$ a non-modular point of $\A$, then $m_p(\A) >m_q(\A)$.

\end{lem}
The following result relates the multiplicity of a modular point $p$ in $\A:f=0$ to the integer $r=mdr(f)$.

\begin{prop}
\label{propA}
If $\A:f=0$ is a supersolvable line arrangement and $p$ is a modular point of $\A$, then either $m_p(\A)=r+1$, or $m_p(\A)=d-r$. In particular, one has 
$$r=\min \{m_p(\A)-1, d-m_p(\A)\}.$$

\end{prop} 

\proof The central projection from $p$ induces a locally trivial fibration with base $B$, equal to $\PP^1$ minus $m_p=m_p(\A)$ points, total space $M(\A)$, the complement of the line arrangement $\A$ in $\PP^2$, and fiber $F$,  obtained from $\PP^1$ by deleting $d-m_p+1$ points. It follows that
\begin{equation} \label{r3}
\chi(M(\A))=\chi(B)\chi(F)=(m_p-2)(d-1-m_p).
\end{equation}
 On the other hand, for a line arrangement $\A$, the global Tjurina number $\tau(\A)$ coincides to the global Milnor number $\mu(\A)$, which is the sum of all local Milnor numbers of the multiple points of $\A$. Indeed, any such singular point is weighted homogeneous, and hence we apply again K. Saito's result, see \cite{KS}. Hence one has
\begin{equation} \label{r4}
\chi(\A)=2-(d-1)(d-2)+\tau(\A).
\end{equation}
Here $\A$ is regarded as a singular plane curve and we use a well known formula, see for instance \cite[Formula (4.5)]{DHA}.
It follows that 
$$\tau(\A)=\chi(\A)+(d-1)(d-2)-2=\chi(\PP^2)-\chi(M(\A))+d^2-3d=$$
$$=(d-1)^2-(m_p-1)(d-1-(m_p-1)).$$
Now $\A$ is free, since it is supersolvable, and hence there is equality in formula \eqref{r1}.
It follows that $r=m_p-1$ or $r=d-m_p$.
\endproof

\begin{ex}
\label{exA}
In the full monomial line arrangement 
$$\A: f=xyz(x^m-y^m)(x^m-z^m)(y^m-z^m)=0,$$
for $m \geq 1$, one has $r=m+1$, see \cite[Example 8.6 (ii)]{DHA}. The modular points are the points of multiplicity $m+2$. Hence in this case they are all of multiplicity $m_p=r+1$. 

\end{ex}

\begin{ex}
\label{exA1}
For two integers $i \leq j$ we define a homogeneous polynomial in $\C[u,v]$ of degree $j-i+1$ by the formula 
\begin{equation}
\label{gij}
 g_{i,j}(u,v)=(u-iv)(u-(i+1)v) \cdots (u-jv).
\end{equation}
Consider the line arrangement $\A : f=0$ of $d=d_1+d_2+1 \geq 3$ lines in $\PP^2$ given by
$$f(x,y,z)= xg_{1,d_1}(x,y)g_{1,d_2}(x,z)=0,$$
for $1 \leq d_1 <d/2$ and $d_2=d-1-d_1$. 
This line arrangement, denoted by $\hat L(d_1+1,d_2+1)$,  was considered in \cite[Example 4.10]{DIM}, \cite[Remark 4.1]{DStExpo}, and  is free with exponents $(d_1,d_2)$. Moreover  it has two modular points, one of multiplicity $m_1=d_1+1=r+1$, the other of multiplicity $m_2=d_2+1=d-r$. Hence both cases in Proposition \ref{propA} can occur.
\end{ex}

\begin{ex}
\label{exA2}
For the monomial line arrangement 
$$\A(m,m,3): f=(x^m-y^m)(x^m-z^m)(y^m-z^m)=0,$$
 for $m \geq 2$, one has $r=m+1$ and moreover the  line arrangement $\A(m,m,3)$ is free with exponents $(m+1,2m-2)$, see \cite[Example 8.6 (i)]{DHA}. 
 There are no modular points,  but just intersection points of multiplicity $3$ and $m$. So the equalities $m_p(\A)=r+1$ and $m_p(\A)=d-r$ can both fail for a free line arrangement which is not supersolvable.

\end{ex}

To refer to certain line arrangements in $\PP^2$, we recall the following notation from \cite{DIM}. We say that a line arrangement $\A$ of $d$ lines is of type $L(d,m)$ if there is a single point of multiplicity $m \geq 3$ and all the the other intersection points of $\A$  are double points. 
We recall also the following result. 

\begin{prop}
\label{propArec}
Let $\A:f=0$ be a  line arrangement of $d$ lines in $\PP^2$. Then one has the following.

\begin{enumerate}
\item $mdr(f)=1$ if and only if $d=3$ and $\A$ is a triangle, or $d \geq 4$ and $\A$ is of type $L(d,d-1)$. Any such arrangement is free.

\item Any arrangement $\A$ of type $L(d,d-2)$ for $d\geq 5$ is nearly free, with $mdr(f)=2$.

\item Any arrangement $\A:f=0$ with $mdr(f)=2$ is either of type $L(d,d-2)$, or of type 
$\hat L(3,m_2)$, or linearly equivalent to the monomial line arrangement $\A(2,2,3)$.

\end{enumerate}

\end{prop}

For the claims (1) and (2) we refer to \cite[Proposition 4.7]{DIM}, and for (3)  we refer to \cite{To} or \cite[Theorem 4.11]{DIM}.

\section{The freeness properties of nearly supersolvable line arrangements}

 The following result is the analog of Lemma \ref{lemA} in this setting.
 \begin{prop}
\label{propA1}
If $\A:f=0$ is a nearly supersolvable line arrangement  and $p$ is a nearly modular point of $\A$, then 
$$m_p(\A)=max \{m_q(\A) : q \in \A \}.$$
\end{prop} 

\proof Consider the supersolvable arrangement $\B=\B(\A,p)$ defined in the Introduction. Let $q$ be a modular point of $\B$, with $q \ne p$. If $q$ is not on the line $L$, then $q$ is a multiple point of $\A$, and moreover, it is a modular point of $\A$. This is impossible, since $\A$ is not supersolvable. Hence $q \in L$, but the line $L$ contains the point $p$, the point $p'$ and some other points of multiplicity $2$ in $\B$. Among all these points, clearly $p$ has the largest multiplicity in $\B$, namely $m_p(\A)+1$. Any other multiple point $q' \notin L$ of $\A$, occurs as a multiple point of $\B$ with the same multiplicity $m_{q'}(\A)=m_{q'}(\B)$.
Lemma \ref{lemA} implies that
$$m_p(\B)=m_p(\A)+1 > m_{q'}(\B)=m_{q'}(\A),$$
and hence $m_p(\A) \geq m_{q'}(\A)$. Since $m_p(\A) \geq 2=m_{p'}(\A)$, the claim is proved.
\endproof

 \begin{prop}
\label{propA2}
Any line arrangement $\A:f=0$ with $r=mdr(f) \leq 2$ is either supersolvable or nearly supersolvable. In particular, if $\A$ consists of $d \leq 5$ lines and if $\A$ is either free or nearly free, then $\A$ is either supersolvable or nearly supersolvable.
\end{prop} 

\proof
Use Proposition \ref{propArec} (3) and note that the arrangements $\hat L(m_1,m_2)$ and $\A(2,2,3)$ are supersolvable, while $L(d,d-2)$ is clearly nearly supersolvable. The last claim follows from the inequality $2r \leq d$.
\endproof

Our interest in this class of line arrangements comes from the following result.

\begin{thm}
\label{thmB}
Let $\A:f=0$ be a nearly supersolvable line arrangement of $d$ lines in $\PP^2$, and let $p$ be a nearly modular point of $\A$. Then either $mdr(f)=d-m_p(\A)$ and then $\A$ is nearly free, or $d=2d_1+1$,  
$mdr(f)=m_p(\A)=d_1$ and then $\A$ is free. In fact, the first case occurs when $2m_p(\A) \geq d$, while the second case occurs when $2m_p(\A) =d-1$.

\end{thm} 

\proof 
We set again $m_p=m_p(\A)$.
The central projection from $p$ induces a locally trivial fibration with base $B$, equal to $\PP^1$ minus $m_p+1$ points, total space $M(\B)$, the complement of the line arrangement $\B$ in $\PP^2$, and fiber $F$,  obtained from $\PP^1$ by deleting $d-m_p+1$ points. It follows that
\begin{equation} \label{r5}
\chi(M(\B))=\chi(B)\chi(F)=(m_p-1)(d-1-m_p).
\end{equation}
 Note that $M(\A)$ is the disjoint union of $M(\B)$ with $L'$, where $L'$ is obtained from $L$ by deleting $d-m_p$ points, and hence
\begin{equation} \label{r6}
\chi(M(\A))=\chi(M(\B))+\chi(L')=(m_p-1)(d-1-m_p)+(2-d+m_p)=
\end{equation}
$$=(m_p-2)(d-1-m_p)+1.$$
It follows as above that 
$$\tau(\A)=\chi(\PP^2)-\chi(M(\A))+d^2-3d=(d-1)^2-(m_p-1)(d-1-(m_p-1))-1.$$
Now we apply \cite[Theorem 1.2]{Dcurves} and we get the following possibilities.

\begin{enumerate}

\item $mdr(f)=d-m_p$. Then the formula \eqref{r2} implies that $\A$ is nearly free.
In particular, in this case $d-m_p \leq d/2$, and hence $m_p\geq d/2$.

\item $mdr(f)=m_p-1$ and $\A$ is free. But the formula \eqref{r2} implies that $\A$ is nearly free, hence we get a contradiction in this case.

\item $m_p \leq mdr(f) \leq d-m_p-1$, and in particular $m_p\leq (d-1)/2$. But then we know that 
$$\tau(\A)\leq (d-1)^2-m_p(d-1-m_p)$$
by using \eqref{r1}. Hence
$$m_p(d-1-m_p) \leq (m_p-1)(d-1-(m_p-1))+1,$$
which implies $m_p \geq (d-1)/2$. Hence this case is possible only when $d=2d_1+1$ is odd, and 
$$r=mdr(f)=m_p=\frac{d-1}{2}=d_1.$$
These equalities imply that
$$\tau(\A)=(d-1)^2-(r-1)(d-r)-1=(d-1)^2-r(d-r-1)$$
and hence in this case $\A$ is a free arrangement.
\end{enumerate}

\endproof
The following direct consequence of Theorem \ref{thmB} is rather surprising, in view of the fact that the multiplicity of a modular point can be arbitrarily small, as shown by Example \ref{exA1}.
\begin{cor}
\label{corS}
Let $\A:f=0$ be a nearly supersolvable line arrangement of $d$ lines in $\PP^2$, and let $p$ be a nearly modular point of $\A$. Then $$m_p(\A) \geq \frac{d-1}{2}.$$

\end{cor}
This result has the following application to the Slope problem, which we recall briefly here following \cite[Subsection (2.2)]{T1}. Let $\K=\R, \C$.
Consider $n\geq 3$ distinct points $P_1,...,P_n \in \K^2$, not all collinear and consider the set of lines $L_{i,j}$ determined by all the pairs of points $P_i,P_j$ for $i<j$. Two such lines $L_{i,j}$ and $L_{i',j'}$ {\it have the same slope} if, when we embed $\K^2$ in the projective space $\PP^2(\K)$, the closures $\overline L_{i,j}$ and $\overline L_{i',j'}$ of the lines $L_{i,j}$ and $L_{i',j'}$ meet the line at infinity
$L=\PP^2(\K) \setminus \K^2$ at the same point $D$. Let $D_1,...,D_w$ be the points on $L$ obtained by taking the intersections with all the 
closures $\overline L_{i,j}$ of the lines $L_{i,j}$ for $1\leq i<j\leq n$.
The Slope problem claims that in these conditions and when $\K=\R$, there are at least $n-1$ slopes, i.e. with our notation, one has 
$$w \geq n-1.$$
It is known that this inequality fails for the case $\K=\C$, see Remark \ref{rkS}.
 If we dualize this setting, as explained in \cite[Subsection (2.2)]{T1}, we replace the points $P_i$ by the dual lines $\ell_i$ in $\PP^2(\K)$, the points $D_j$ by the lines $\delta_j$ and the line at infinity $L$ becomes a point $P_L$.
The line arrangement $\A=\{\ell_1,...,\ell_n, \delta_1, ..., \delta_w\}$ is supersolvable, with $P_L$ a modular point of multiplicity $w$, since all the lines $\delta_j$ pass through $P_L$. 
\begin{thm}
\label{thmSlope2}
With the above notation, assume that one of the points $D_k$, say the point $D_1$, is obtained as an intersection $L \cap \overline L_{i,j}$ for a unique pair $i<j$. Then one has the stronger inequality $$w\geq n,$$
valid in both cases $\K=\R, \C$.
\end{thm}
\proof
The proof is just a variation of the proof of  \cite[Proposition 2.7]{T1}, in which the supersolvable arrangements are replaced by nearly supersolvable arrangements.
It is enough to note that the line arrangement $\B$ obtained from the above line arrangement $\A=\{\ell_1,...,\ell_n, \delta_1, ..., \delta_w\}$
by deleting the line $\delta_1$ is nearly supersolvable with $P_L$ as its nearly modular point and the unique intersection point not on a line through $P_L$ in $\B$ being the intersection $\ell_i \cap \ell_j \in \delta_1$. Indeed, the lines $\ell_i , \ell_j, \delta_1$ meet since the dual points $P_i,P_j,D_1$ are collinear.
We conclude by applying Corollary \ref{corS}.
\endproof

\begin{rk}
\label{rkS} When $\A:f=0$ is a  supersolvable line arrangement of $d$ lines in $\PP^2$, and $p$ is a modular point of $\A$, then the inequality 
$$m_p(\A) \geq \frac{d-1}{2}$$
can fail. To see this, consider the real supersolvable line arrangement from Example \ref{exA1} with $d_1<d_2-2$ and $p$ the modular point of multiplicity $m_1=d_1+1$. On the other hand, if we set
$$m(\A)=max \{m_q(\A) : q \in \A \},$$
then, for a real supersolvable  arrangement  $\A$ of $d$ lines, one has the inequality 
$$m(\A) \geq \frac{d-1}{2},$$
see \cite[Proposition 2.7]{T1}, where it is shown that this inequality is equivalent to (a positive answer to) the Slope problem.
Example \ref{exA} shows that the inequality $m(\A) \geq \frac{d-1}{2}$
fails for the (complex supersolvable) full monomial line arrangement 
for $d=|\A|=3m \geq 9$.
\end{rk}

\begin{rk}
\label{rkNF2}
If $\A:f=0$ is a nearly supersolvable line arrangement of $d$ lines in $\PP^2$ and let $p$ be a nearly modular point of $\A$. Assume that $2m_p(\A) \geq d+1$ and hence the line arrangement $\A$ is nearly free. As mentioned above, $\A$ is obtained by deletion of one line $L$ from the free arrangement $\B$ of $d+1$ lines, with exponents $d_1=d-m_p(\A)$ and $d_2=m_p(\A)$. Note that the line $L$ contains only two multiple points of $\B$, namely $p$ with multiplicity $m_p(\A)+1$ and $p'$ with multiplicity $3$, i.e. a triple point. A distinct construction of a nearly free line arrangement $\A$ from a free arrangement $\B$ with exponents $(d_1,d_2)$ by deleting one line $L'$ is presented in
 \cite[Proposition 3.1]{MaVa}. In their construction, the line $L'$ should contain $t=d_2$ triple points of $\B$, and hence our Theorem \ref{thmB} does not cover a special case of their construction. Note however that in both cases, the arrangement $\A$ inherits the exponents of $\B$.
\end{rk}

\begin{ex}
\label{exRef1}
Let $\tilde A(m_1,m_2):f=0 $ be the line arrangement obtained by taking the union of two pencils of lines in $\PP^2$, one containing $m_1\geq 2 $ lines, the other containing $m_2 \geq m_1$ lines, in general position to each other.
Then it is easy to see that $d=m_1+m_2$, $mdr(f)=m_1$ and 
$$\tau(\tilde A(m_1,m_2))=(d-1)^2-m_1(d-m_1-1)-(m_1-1),$$
see \cite[Proposition 4.9]{DIM}.
Choose $p$ to be the base point of the second pencil, hence a point of multiplicity $m_2$, and $p'$ to be the base point of the first  pencil, hence a point of multiplicity $m_1$. Then for $m_1=2$ we get a nearly supersolvable arrangement, with $p$ a nearly modular point, and the above formula combined with the equality \eqref{r2} implies that this arrangement is nearly free. Note that for $m_1>2$ the arrangement $\tilde A(m_1,m_2)$ is neither free, nor nearly free. This fact explains why in the definition of a nearly modular point we have considered only points $p'$ of multiplicity 2.
\end{ex}

\begin{ex}
\label{exB}

Consider the line arrangement $\A : f=0$ of $d=d_1+d_2 \geq 4$ lines in $\PP^2$ given by
$$f(x,y,z)= x(y-z)g_{1,d_1-1}(x,y)g_{2,d_2}(x,z)=0$$
for $2 \leq d_1 \leq d_2$, with $g_{i,j}$ as in Example \ref{exA1}. The point $p_1=(0:0:1)$ has multiplicity $d_1$ (the number of factors in $f$ involving only the variables $x$ and $y$), and the point $p_2=(0:1:0)$ has multiplicity $d_2$.
The line arrangement $\A$ is nearly free with exponents $(d_1,d_2)$, see \cite[Example 4.14]{DIM} or \cite[Proposition 4.2]{DStExpo}. We describe now a minimal set of generators for the graded $S$-module $AR(f)$. Define new homogeneous polynomials $A(u,v)$ of degree $d_1-1$, $B(u,v)$ of degree $d_2-1$, $C(u,v)$ of degree $d_1-2$ and $D(u,v)$ of degree $d_2-d_1+1$ by the following relation
$$A(u,v)=g_{1,d_1-1}(u,v)=(u-v)C(u,v) \text{ and } B(u,v)=g_{2,d_2}(u,v)=C(u,v)D(u,v).$$
Define now the syzygies $\rho_i=(a_i,b_i,c_i) \in AR(f)$ for $i=1,2,3$ as follows.
The syzygy $\rho_1$ has degree $d_1$ and is given by
$$a_1=x(A(x,y)+(y-z)A_y(x,y)),$$
$$b_1= y(A(x,y)+(y-z)A_y(x,y))-d(y-z)A(x,y),$$
$$c_1=z(A(x,y)+(y-z)A_y(x,y).$$
The syzygy $\rho_2$ has degree $d_2$ and is given by
$$a_2=x(-B(x,z)+(y-z)B_z(x,z)), $$
$$b_2=y(-B(x,z)+(y-z)B_z(x,z)), $$
$$c_2=z(-B(x,z)+(y-z)B_z(x,z))- d(y-z)B(x,z).$$
The syzygy $\rho_3$ has again degree $d_2$ and is given by
$$a_3=x(x-2y+z)[D(x,z)C(x,y,z)+C_z(x,z)D(x,z)+C(x,z)D_z(x,z)] +$$
$$+x(x-y)C_y(x,y)D(x,z),$$
$$b_3=y(x-2y+z)[D(x,z)C(x,y,z)+C_z(x,z)D(x,z)+C(x,z)D_z(x,z)] +$$
$$+d(y-x)C(x,y)D(x,z)+ y(x-y)C_y(x,y)D(x,z),$$
and
$$c_3=z(x-2y+z)[D(x,z)C(x,y,z)+C_z(x,z)D(x,z)+C(x,z)D_z(x,z)]-$$
$$-d(x-2y+z)C(x,z)D(x,z)+z(x-y)C_y(x,y)D(x,z).$$
In the above formulas a subscript indicate a partial derivatives, for instance $A_y(x,y)$ is the partial derivatives of $A(x,y)$ with respect to $y$. More we use the notation
$$C(x,y,z)=\frac{C(x,y)-C(x,z)}{y-z}.$$
The first (resp. second) syzygy $\rho_1$ (resp. $\rho_2$) is obtained using the 
multiple point $p_2$ (resp. $p_1$) and the general recipe presented in \cite[Section (2.2)]{Dcurves}. The third syzygy $\rho_3$ is obtained by dividing the vector in $S^3$
$$D(x,z) \rho_1 +(x-2y+z)\rho_2$$
by $y-z$. Hence we get the following generator of the relations among $\rho_1,\rho_2$ and $\rho_3$:
\begin{equation} \label{Syz2}
D(x,z) \rho_1 +(x-2y+z)\rho_2-(y-z)\rho_3=0.
\end{equation}
Note that $\A$ is  nearly supersolvable, with a nearly modular point given by $p=(0:1:0)$ with multiplicity $m_p=d_2$. The only node not connected by a line to $p$ is the point $p'=(1:1:1)$.
If $d_1 <d_2$ or if $d_1=d_2$ and the exact sequence in Proposition \ref{propNF} or in Remark \ref{rkNF} is constructed using $\rho_1, \rho_2, \rho_3$ above, then $p'=P(\A)$ is the jumping point of the arrangement, i.e. the solution of the equations
$$x-2y+z=y-z=0$$
coming from \eqref{Syz2}, see also \cite[Proposition 2.7]{MaVa}. 
\end{ex}

\begin{rk}
\label{rkNF3}
The nearly free arrangement $\A:f=0$ in Exemple \ref{exB} is obtained from the free arrangement $\C:g(x,y,z)= xg_{1,d_1-1}(x,y)g_{2,d_2}(x,z)=0$ with exponents $(d_1'=d_1-1,d'_2=d_2-1)$, of type $\hat L(d_1,d_2)$ as discussed in Example \ref{exA1}, by adding the line $L:y-z=0$. Note that this line $L$ contains $d_1-2=d'_1-1$ triple points, namely the points $(k:1:1)$ for $k=2,3,...,d_1-1$, and hence the arrangement $\A$ can be regarded as the result of the construction described in 
\cite[Proposition 3.3]{MaVa}. In particular, \cite[Proposition 3.5]{MaVa} implies that the jumping point $P(\A)$ should belong to the line $L$, a result less precise than what we have shown above by explicit computation, namely that $p'=P(\A)$. Moreover,  in the case $d_1=d_2=2$  the software SINGULAR gives different generators $\rho_1',\rho_2',\rho_3'$ for $AR(f)$, namely
$$\rho_1'= (4xy-5xz, 4xy-4xz-yz, 4xy-4xz-4yz+3z^2 ),$$
$$\rho_2'= (x^2-2xy+xz,-3xy+2y^2+4xz-3yz, xz-2yz+z^2)$$
and
$$
\rho_3'= (12x^2+4xy-23xz,12x^2+4xy-28xz+5yz, 12x^2+4xy-28xz-4yz+9z^2)$$
One checks the following relation
$$ ( -3x-y+4z)\rho_1'+
  (-3z )\rho_2'+
  ( y-z)\rho_3'=0$$
 Depending which of the syzygies $\rho_1'$, $\rho_2'$ and $\rho_3'$ are chosen as the first syzygy $\rho_1$ in the exact sequence from Proposition \ref{propNF}, we get the following jumping points:
 
\noindent (i) $z=y-z=0$, hence $P(\A)=(1:0:0)$, when $\rho_1=\rho_1'$,
 
\noindent (ii) $-3x-y+4z=y-z=0$, hence $P(\A)=(1:1:1)$, when $\rho_1=\rho_2'$, 

\noindent (iii) $-3x-y+4z=z=0$, hence $P(\A)=(1:-3:0)$, when $\rho_1=\rho_3'$.

Hence when $d_1=d_2$ the choice of the jumping point $P(\A)$ is not unique.

\end{rk}

\begin{ex}
\label{exC}
Consider the line arrangement $\A : f=0$ of $d=2d_1+1 \geq 5$ lines in $\PP^2$ given by
$$f(x,y,z)= x(y-z)(y+x-z)g_{1,d_1-1}(y,x)g_{2,d_1}(z,x)=0$$
for $2 \leq d_1$, with $g_{i,j}$ as in Example \ref{exA1}. Then the line arrangement $\A$ is  free with exponents $(d_1,d_1)$ by Theorem \ref{thmB}. Indeed, $\A$ is  nearly supersolvable, with a nearly modular point given by $p=(0:1:0)$ with multiplicity $m_p=d_1$. The only node not connected by a line to $p$ is the point $p'=(1:1:1)$. It follows that 
this example corresponds to the case (3) in the  proof of Theorem \ref{thmB}, since $d_1=mdr(f)=(d-1)/2$.
Note that this arrangement has a second nearly modular point at $p_0=(0:0:1)$, with the corresponding node at $p_0'=(1:d_1:d_1)$. This shows that the nearly modular point is not necessarily unique when $d_1=d_2$.
\end{ex}

\begin{cor}
\label{corA}
Any nearly supersolvable free (resp. nearly free) line arrangement $\A$ satistfies Terao's Conjecture, namely
if another line arrangement $\B$ has the same intersection lattice as $\A$, then $\B$ is also free (resp. nearly free) with the same exponents as the line arrangement $\A$.
\end{cor}

\proof
In the above statement, the intersection lattices refers in fact to the intersection lattices of the corresponding central plane arrangements in $\C^3$. It is clear that nearly supersolvability is a combinatorial property, and hence $\B$ is also nearly supersolvable. If $p$ (resp. $q$) denotes  a nearly modular point for $\A$ (resp. the corresponding  nearly modular point for $\B$), then clearly $m_p(\A)=m_q(\B)$. The claim follows then from Theorem \ref{thmB}.

\endproof

\begin{rk}
\label{rkNF4}
If $\A:f=0$ is a  supersolvable line arrangement of $d$ lines in $\PP^2$, then it is known that the complement $M(\A)$ is a $K(\pi,1)$-space, see \cite[Theorems 4.15 and 4.16]{DHA}.
When $\A:f=0$ is the nearly supersolvable line arrangement from Example \ref{exB}, then
\cite[Proposition 5.3]{DStExpo} shows that the complement $M(\A)$ is not a $K(\pi,1)$-space,
at least when $d_1=2$ or for the pair $(d_1,d_2)=(3,3)$.
\end{rk}

\section{Splitting types and jumping lines for the bundle of logarithmic vector fields} 

Let $E_{\A}$ be the locally free sheaf on $X=\PP^2$ defined by
\begin{equation} \label{equa1} 
E_{\A}=T\langle \A \rangle (-1),
\end{equation}
where $T\langle \A \rangle $ is the sheaf of logarithmic vector fields along $\A$ as considered for instance in \cite{AD, DS14, DStJump}.
For a line $L$ in $X$, the pair of integers $( d_1^L, d_2^L)$, with $ d_1^L \leq d_2^L$, such that $ E_{\A}|_L \simeq \OO_L(-d_1^L) \oplus \OO_L(-d_2^L)$ is called the splitting type of $E_{\A}$ along $L$, see for instance \cite{FV1,OSS}. For a generic line $L_0$, the corresponding splitting type $( d_1^{L_0}, d_2^{L_0})$ is constant.

 Note that  $\A$ is free with exponents $d_1\leq d_2$  if and only if $E_{\A}=\OO_X(-d_1) \oplus \OO_X(-d_2)$, and hence the splitting type is $(d_1,d_2)$ for any line $L$.
One has the following result, see \cite[Lemma 3.6]{DMO}, or apply Lemma \ref{lemA} and Proposition \ref{propA} above.

\begin{cor}
\label{corB}
Let $\A$ be a supersolvable line arrangement, with $d=|\A|$ and $m=m(\A)$ the maximal multiplicity of an intersection point of $\A$. Then the (unordered) splitting type of $\A$ along any line $L$ is $(m-1,d-m)$.
\end{cor}

When  $\A$ is nearly free with exponents $d_1\leq d_2$, then the generic (unordered) splitting type is $(d_1,d_2-1)$ and an unordered splitting type is $(d_1-1,d_2)$ for a jumping line $L$, see \cite[Corollary 3.4]{AD}. This implies the following via Theorem \ref{thmB}.

\begin{cor}
\label{corC}
Let $\A:f=0$ be a nearly supersolvable line arrangement, with $d=|\A|$ and $p$  a nearly modular intersection point  of $\A$. Then the only possible cases are the following.
\begin{enumerate}

\item $2m_p(\A)<d$. Then 
$d=2m_p(\A)+1$ is odd and  $\A$ is free with exponents $(d_1,d_1)$ with $d_1=m_p(\A)$, the generic splitting type is $(d_1,d_1)$ and there are no jumping lines.

\item $2m_p(\A) = d$. Then $\A$ is nearly free with exponents $(d_1,d_1)$ with
$d_1=m_p(\A)$,
the generic splitting type is $(d_1-1,d_1)$ and there are no jumping lines.

\item $2m_p(\A) > d$. Then $\A$ is nearly free with exponents $(d_1,d_2)$ with
$d_1=d-m_p(\A)$, $d_2=m_p(\A)$ and the generic splitting type is $(d_1,d_2-1)$.
A line $L$ is a jumping line if and only if it passes through the jumping point $P(\A)$ of the nearly free arrangement $\A$, and the corresponding splitting type is 
 $(d_1-1,d_2)$.

\end{enumerate} 
In particular, for a nearly supersolvable line arrangement the generic splitting type of $T\langle \A \rangle $ is determined by the combinatorics.

\end{cor}

\proof
The only claim that needs justification is the last one, which follows from \cite[Proposition 2.4]{MaVa}. Note that in this case $d_1<d_2$ and hence the jumping point $P(\A)$ is uniquely defined by the line arrangement $\A$ as we noticed in Proposition \ref{propNF} and Remark \ref{rkNF3}.

\endproof
One can use this result and \cite[Theorem 1.2]{CHMN} to show that a finite set of points $Z$ in $\PP^2$ whose dual line arrangement $\A_Z$ is nearly supersolvable never admits an unexpective curve. We refer to \cite{CHMN, DMO} for more on this subject, see in particular  \cite[Theorem 3.7]{DMO}.

\begin{rk}
\label{rkNF5}
It is a major open question whether the generic splitting type of $T\langle \A \rangle $ is determined by combinatorics for any line arrangement, see \cite[Question 7.12]{CHMN}.
The nearly supersolvable line arrangements form a class where this question has a positive answer. For the moment there is no combinatorial description for the larger class of nearly free line arrangements, hence if $\A$ is nearly free and $\A'$ has the same combinatorics as $\A$,
we do not know whether $\A'$ is also nearly free. When this is the case, then $\A$ and $\A'$ have the same exponents and hence the same generic splitting type for $T\langle \A \rangle $
and for $T\langle \A '\rangle $.
\end{rk}

\end{document}